\documentclass{article}

\usepackage{amssymb,latexsym,amsmath}

\usepackage{graphicx}

\usepackage{graphicx,amssymb,lineno}

\begin{document}

\newcommand{\bfi}{\bfseries\itshape}

\makeatletter

\@addtoreset{figure}{section}

\def\thefigure{\thesection.\@arabic\c@figure}

\def\fps@figure{h, t}

\@addtoreset{table}{bsection}

\def\thetable{\thesection.\@arabic\c@table}

\def\fps@table{h, t}

\@addtoreset{equation}{section}

\def\theequation{\thesubsection.\arabic{equation}}

\makeatother

\newtheorem{thm}{Theorem}[section]

\newtheorem{prop}[thm]{Proposition}

\newtheorem{lema}[thm]{Lemma}

\newtheorem{cor}[thm]{Corollary}

\newtheorem{defi}[thm]{Definition}

\newtheorem{rk}[thm]{Remark}

\newtheorem{exempl}{Example}[section]

\newenvironment{exemplu}{\begin{exempl}  \em}{\hfill $\surd$

\end{exempl}}

\newcommand{\comment}[1]{\par\noindent{\raggedright\texttt{#1}

\par\marginpar{\textsc{Comment}}}}

\newcommand{\todo}[1]{\vspace{5 mm}\par \noindent \marginpar{\textsc{ToDo}}\framebox{\begin{minipage}[c]{0.95 \textwidth}

\tt #1 \end{minipage}}\vspace{5 mm}\par}

\newcommand{\ea}{\mbox{{\bf a}}}

\newcommand{\eu}{\mbox{{\bf u}}}

\newcommand{\ueu}{\underline{\eu}}

\newcommand{\ueo}{\overline{u}}

\newcommand{\oeu}{\overline{\eu}}

\newcommand{\ew}{\mbox{{\bf w}}}

\newcommand{\ef}{\mbox{{\bf f}}}

\newcommand{\eF}{\mbox{{\bf F}}}

\newcommand{\eC}{\mbox{{\bf C}}}

\newcommand{\en}{\mbox{{\bf n}}}

\newcommand{\eT}{\mbox{{\bf T}}}

\newcommand{\eL}{\mbox{{\bf L}}}

\newcommand{\eR}{\mbox{{\bf R}}}

\newcommand{\eV}{\mbox{{\bf V}}}

\newcommand{\eU}{\mbox{{\bf U}}}

\newcommand{\ev}{\mbox{{\bf v}}}

\newcommand{\eve}{\mbox{{\bf e}}}

\newcommand{\uev}{\underline{\ev}}

\newcommand{\eY}{\mbox{{\bf Y}}}

\newcommand{\eK}{\mbox{{\bf K}}}

\newcommand{\eP}{\mbox{{\bf P}}}

\newcommand{\eS}{\mbox{{\bf S}}}

\newcommand{\eJ}{\mbox{{\bf J}}}

\newcommand{\eB}{\mbox{{\bf B}}}

\newcommand{\eH}{\mbox{{\bf H}}}

\newcommand{\leb}{\mathcal{ L}^{n}}

\newcommand{\eI}{\mathcal{ I}}

\newcommand{\eE}{\mathcal{ E}}

\newcommand{\hen}{\mathcal{H}^{n-1}}

\newcommand{\eBV}{\mbox{{\bf BV}}}

\newcommand{\eA}{\mbox{{\bf A}}}

\newcommand{\eSBV}{\mbox{{\bf SBV}}}

\newcommand{\eBD}{\mbox{{\bf BD}}}

\newcommand{\eSBD}{\mbox{{\bf SBD}}}

\newcommand{\ecs}{\mbox{{\bf X}}}

\newcommand{\eg}{\mbox{{\bf g}}}

\newcommand{\paromega}{\partial \Omega}

\newcommand{\gau}{\Gamma_{u}}

\newcommand{\gaf}{\Gamma_{f}}

\newcommand{\sig}{{\bf \sigma}}

\newcommand{\gac}{\Gamma_{\mbox{{\bf c}}}}

\newcommand{\deu}{\dot{\eu}}

\newcommand{\dueu}{\underline{\deu}}

\newcommand{\dev}{\dot{\ev}}

\newcommand{\duev}{\underline{\dev}}

\newcommand{\weak}{\stackrel{w}{\approx}}

\newcommand{\mild}{\stackrel{m}{\approx}}

\newcommand{\strong}{\stackrel{s}{\approx}}

\newcommand{\weakdown}{\rightharpoondown}

\newcommand{\opg}{\stackrel{\mathfrak{g}}{\cdot}}

\newcommand{\opunu}{\stackrel{1}{\cdot}}
\newcommand{\opdoi}{\stackrel{2}{\cdot}}

\newcommand{\opn}{\stackrel{\mathfrak{n}}{\cdot}}

\newcommand{\tr}{\ \mbox{tr}}

\newcommand{\Ad}{\ \mbox{Ad}}

\newcommand{\ad}{\ \mbox{ad}}

\renewcommand{\contentsname}{ }

\title{Blurred maximal cyclically monotone sets and 
 bipotentials}

\author{ 
Marius Buliga\footnote{"Simion Stoilow" Institute of Mathematics of the Romanian Academy,
 PO BOX 1-764,014700 Bucharest, Romania, e-mail: Marius.Buliga@imar.ro }, 
G\'ery de Saxc\'e\footnote{Laboratoire de 
  M\'ecanique de Lille, UMR CNRS 8107, Universit\'e des Sciences et 
  Technologies de Lille,
 B\^atiment Boussinesq, Cit\'e Scientifique, 59655 Villeneuve d'Ascq cedex, 
 France, e-mail: gery.desaxce@univ-lille1.fr},
Claude  Vall\'ee\footnote{Laboratoire de 
M\'ecanique des Solides, UMR 6610, UFR SFA-SP2MI, bd M. et P. Curie, 
t\'el\'eport 2, BP 30179, 86962 Futuroscope-Chasseneuil cedex, 
France, e-mail: vallee@lms.univ-poitiers.fr}
}

\date{ }

\maketitle

\begin{abstract}

Let $X$ be a reflexive Banach space and $Y$ its dual. 
In this paper we find necessary and sufficient conditions for the existence of a
bipotential for a blurred maximal cyclically monotone set. Equivalently, we 
find a necessary and sufficient condition on $\phi \in \Gamma_{0}(X)$ for 
that the differential inclusion 
 $\displaystyle y \in \bar{B}(\varepsilon) + \partial \phi(x)$ can be put 
in the form $y \in \partial b(\cdot, y)(x)$, with $b$ a bipotential. 

\end{abstract}

{\bf MSC-class:} 49J53; 49J52; 26B25


\section{Introduction}

Let $\phi \in \Gamma_{0}(X)$, where $X$ is a reflexive Banach space, with 
dual $Y$, and $\varepsilon > 0$. 
We want to  describe the set $M(\phi, \varepsilon)$  of
solutions $(x,y) \in X \times Y$ of the following problem:  
\begin{equation}
 \mbox{ there 
is } a \in Y \, , \, \|a\| \leq \varepsilon \mbox{  such that } 
 y + a \in \partial \phi(x) 
\label{p}
\end{equation}

 Our main result are proposition \ref{pnewc} and  theorem \ref{maithm}, which
 imply the following (for the notion of bipotential see definition \ref{def2}). 

\begin{cor}
If  for any $y \in Y$ the set $\displaystyle 
\bigcup_{\| \bar{y}-y\|\leq \varepsilon} \partial 
\phi^{*}(\bar{y})$  is convex then there exists a bipotential  
$b: X \times Y \rightarrow \bar{\mathbb{R}}$ such that $M(\phi, \varepsilon) =
M(b)$. Therefore, there is a   lower semicontinuous and convex
in each variable function $b$, such that   any of the following relations 
\begin{enumerate}
\item[(a)] $y \in \partial b(\cdot, y) (x)$, 
\item[(b)] $x \in \partial b(x, \cdot)(y)$, 
\item[(c)] $b(x,y) = \langle x, y \rangle$, 
\end{enumerate} 
is equivalent with $(x,y) \in M(\phi, \varepsilon)$. 
\label{coro1}
\end{cor}

\paragraph{Proof.}
Let $M = Graph(\partial \phi)$ and $A = \left\{ 0 \right\} \times
\bar{B}_{Y}(\varepsilon)$. Then $M(\phi, \varepsilon) =  M + A$. By proposition
\ref{pnewc} in the 
hypothesis of the corollary (which will appear further as relation (\ref{newc}))
the set $M + A$ is a BB-graph, therefore by theorem \ref{thm1} there exists a bipotential $b$ (definition
\ref{def2}) such that 
$$\displaystyle M+A = M(b) = \left\{ (x,y) \in X \times Y \mbox{ : } b(x,y) = \langle 
x, y \rangle \right\}$$
therefore (\ref{p}) is equivalent also with $y \in \partial b(\cdot, y) (x)$, or 
with $x \in \partial b(x, \cdot)(y)$. $\quad \square$

\begin{cor}
If  the function $b$ defined by  
 \begin{equation}
b(x,y) = \phi(x) + \inf_{\|a\| \leq \varepsilon} \left[ \phi^{*}(y-a) +
\langle x, a \rangle \right]
\label{bemi}
\end{equation}
is convex in the first variable,  then $b$ is a bipotential such that 
$M(b) = M(\phi, \varepsilon)$. 
\label{coro2}
\end{cor}

\paragraph{Proof.}  If $b$ defined by (\ref{bemi}) 
is convex then the function $\displaystyle f: \bar{B}(\varepsilon) \times X \times Y
\rightarrow \bar{\mathbb{R}}$ defined by 
$$f(a, x, y) = \phi(x) + \phi^{*}(y - a) + \langle x, a \rangle$$ 
is implicitly convex in the first two arguments, in the sense of definition 
\ref{defimpl}. Indeed, suppose that  for any $y \in Y$, for any 
 $\displaystyle x_{1}, x_{2} \in X$ and any $\alpha, \beta \in [0,1]$, $\alpha +
 \beta = 1$, such that $\displaystyle b(x_{i}, y) < + \infty$, $i = 1,2$, we
 have 
 $$b(\alpha x_{1} + \beta x_{2}, y) \leq \alpha b(x_{1}, y) + \beta b(x_{2},
 y)$$ 
 The function $\displaystyle f(\cdot, \alpha x_{1} + \beta x_{2}, y)$ is 
 weakly  lower semicontinuous and $\displaystyle \bar{B}(\varepsilon)$ is
 weakly  compact, therefore there exists $\displaystyle \bar{a} \in
 \bar{B}(\varepsilon)$ such that 
 $$f(\bar{a}, \alpha x_{1} + \beta x_{2}, y) = b(\alpha x_{1} + \beta x_{2},
 y)$$
 But then for any $\displaystyle a_{1}, a_{2} \in \bar{B}(\varepsilon)$ we have 
 $$f(\bar{a}, \alpha x_{1} + \beta x_{2}, y) \leq \alpha f(a_{1}, x_{1}, y) + 
 \beta f(a_{2}, x_{2}, y)$$ 
 which proves the previous claim. We use one implication from  theorem
 \ref{maithm} 
 in order to finish the proof. $\quad \square$

\paragraph{Motivation of the problem (\ref{p}).} Recently \cite{bipo4} we have
found a new application of the bipotential method to blurred constitutive laws. 
This application led us to the mathematical problem (\ref{p}) of describing 
blurred maximal cyclically monotone sets by bipotentials. 

The notion of bipotential (definition \ref{def2}) has been introduced in 
\cite{saxfeng}, in order to formulate a large family of non associated 
constitutive laws in terms of convex analysis. The basic idea is explained 
further in few words. In Mechanics the associate constitutive laws are simply 
relations $y \in \partial \phi (x)$, with $\phi: X \rightarrow \mathbb{R} \cup 
\left\{ + \infty \right\}$ a convex and lower semicontinuous function. By
Fenchel inequality such a relation is equivalent with $\displaystyle 
\phi(x) + \phi^{*}(y) = \langle x, y \rangle$, where $\phi^{*}$ is the Fenchel
conjugate of $\phi$. It has been noticed that often in the mathematical 
study of problems related to associated constitutive laws enters not the 
function $\phi$, but the expression 
$$ b(x,y) = \phi(x) + \phi^{*}(y)$$
which we call "separable bipotential". The idea is then to use as 
a basic notion the one of bipotential $b: X \times Y \rightarrow \mathbb{R} \cup 
\left\{ + \infty \right\}$, which is convex and lsc in each argument and
satisfies a generalization of the Fenchel inequality. To non associated
constitutive laws thus corresponds bipotentials which are not separable.

Examples of such laws which can 
be studied with the help of bipotentials are:
non-associated Dr\"ucker-Prager \cite{sax boush KIELCE 93}  and 
Cam-Clay models \cite{sax BOSTON 95} in soil mechanics, 
cyclic Plasticity (\cite{sax CRAS 92},\cite{bodo sax EJM 01}) 
and Viscoplasticity \cite{hjiaj bodo CRAS 00} of metals with non linear 
kinematical hardening rule, Lemaitre's damage law \cite{bodo}, the coaxial 
laws (\cite{sax boussh 2},\cite{vall leri CONST 05}), the Coulomb's friction law \cite{saxfeng}, 
\cite{sax CRAS 92}, 
\cite{boush chaa IJMS 02}, \cite{feng hjiaj CM 06}, \cite{fort hjiaj CG 02}, 
\cite{hjiaj feng IJNME 04}, \cite{sax boush KIELCE 93}, \cite{sax feng IJMCM 98}, \cite{laborde}.
A complete survey can be found in \cite{sax boussh 2}.  

Later we started in \cite{bipo1} \cite{bipo2} \cite{bipo3} a mathematical 
study of bipotentials and their relation with convex analysis. This paper is 
another contribution along this subject. 

Blurred constitutive laws (in mathematical terms blurred graphs of multivalued
operators) appear in many practical situations, due either to experimental 
or numerical indeterminacies \cite{LAD7} \cite{LAD8} \cite{ER3} \cite{ER4}
\cite{ER5}. It is then interesting to take the indeterminacy
into account and to associate, for example, to a differential relation like 
$y \in \partial \phi(x)$ another differential relation $y \in \partial 
b(\cdot, y)(x)$, where $b$ is a bipotential constructed from $\phi$ and the 
indeterminacy. 

We achieve this in the paper,  by using lagrangian convex covers 
introduced in \cite{bipo1}.  For future study is left the more general 
problem of the existence and construction of a bipotential for a blurred 
graph of a multivalued operator which can be expressed by a bipotential. 
Such operators may be monotone, but not cyclically monotone, or even non
monotone. For example in \cite{bipo4} we considered the operator associated 
to Coulomb friction law, which is not even monotone, and we were able to 
construct a differential inclusion for the Coulomb friction law with
indeterminacy.

\paragraph{Aknowledgements.} The first two authors acknowledge the support 
from the European Associated Laboratory "Math Mode" associating the 
Laboratoire de Math\'ematiques de l'Universit\'e Paris-Sud (UMR 8628) and the 
"Simion Stoilow" Institute of Mathematics of the Romanian Academy. We express
our thanks to the anonymous referee for various suggestions leading to
the improvement of the paper. 

\section{Bipotentials and syncs}

$X$ and $Y$ are topological, locally convex, real vector spaces of dual 
variables $x \in X$ and $y \in Y$, with the duality product 
$\langle \cdot , \cdot \rangle : X \times Y \rightarrow \mathbb{R}$. 
We shall suppose that $X, Y$ have topologies compatible with the duality 
product, that is: any  continuous linear functional on $X$ (resp. $Y$) 
has the form $x \mapsto \langle x,y\rangle$, for some $y \in Y$ (resp. 
$y \mapsto \langle x,y\rangle$, for some  $x \in X$). 
We use the notations and conventions of Moreau \cite{moreau}: 
\begin{enumerate} 
\item[-] $\displaystyle \bar{\mathbb{R}} = \mathbb{R}\cup 
\left\{ +\infty \right\}$; by convention $a + (+ \infty) = + \infty$ for any 
$\displaystyle a \in \bar{\mathbb{R}}$ and $a \left(+\infty \right) = + \infty$ 
for any $a \geq 0$; 
\item[-] the domain of a function  $\displaystyle \phi: X \rightarrow
\bar{\mathbb{R}}$ is $dom \, \phi = \left\{ x \in X \mbox{ : } \phi(x) \in
\mathbb{R} \right\}$; 
\item[-] $\displaystyle \Gamma_{0}(X) = \left\{ \phi: X \rightarrow
\bar{\mathbb{R}} 
\mbox{ : } \phi \mbox{ is lsc and } dom \phi \not = \emptyset \right\}$; 
\item[-] for any convex and closed set $A \subset X$, its  indicator function,  
$\displaystyle \chi_{A}$, is defined by 
$$\chi_{A} (x) = \left\{ \begin{array}{ll}
0 & \mbox{ if } x \in A \\ 
+\infty & \mbox{ otherwise } 
\end{array} \right. $$
\item[-] the subdifferential of a function $\displaystyle \phi: X \rightarrow
\bar{\mathbb{R}}$ at a point $x \in X$ is the (possibly empty) set: 
$$\partial \phi(x) = \left\{ u \in Y \mid \forall z \in X  \  \langle z-x, u \rangle \leq \phi(z) - \phi(x) \right\} \  .$$ 
\end{enumerate}

A non empty set $M \subset X \times Y$ is cyclically monotone if for any natural number 
$n \geq 1$ and for any collection $\displaystyle \left\{ (x_{k}, y_{k}) \in M 
\mbox{ : } k = 0, 1, ..., n\right\}$ we have 
$$\langle x_{n} - x_{0} , y_{n}\rangle  + \sum_{1}^{n} \langle x_{k-1} - x_{k}, y_{k-1}
\rangle \, \geq 0 \quad .$$

\begin{defi} A {\bf bipotential} is a function $b: X \times Y \rightarrow
 \bar{\mathbb{R}}$, with the properties: 
\begin{enumerate}
\item[(a)] for any $x \in X$, if $dom \, b(x, \cdot) \not = \emptyset$ then 
$\displaystyle b(x, \cdot) \in \Gamma_{0}(X)$;  for any $y \in Y$, if $dom \,
b(\cdot, y) \not = \emptyset$ then 
$\displaystyle b(\cdot, y) \in \Gamma_{0}(Y)$; 
 \item[(b)] for any $x \in X , y\in Y$ we have $\displaystyle b(x,y) \geq \langle x, y \rangle$; 
\item[(c)]  for any $(x,y) \in X \times Y$ we have the equivalences: 
\begin{equation}
y \in \partial b(\cdot , y)(x) \ \Longleftrightarrow \ x \in \partial b(x, \cdot)(y)  \ \Longleftrightarrow \ b(x,y) = 
\langle x , y \rangle \ .
\label{equiva}
\end{equation}
\end{enumerate}
The {\bf graph} of $b$ is 
\begin{equation}
M(b) \ = \ \left\{ (x,y) \in X \times Y \ \mid \ b(x,y) = \langle x, y \rangle \right\} \  .
\label{mb}
\end{equation}
\label{def2}
\end{defi}

If $X$ is a Banach space,  $Y = X^{*}$ and $M \subset X \times X^{*}$ 
then, by an important result of Rockafellar 
Theorem B \cite{rocka}, is cyclically maximal monotone  if and only if 
$$M \, = \, Graph(\partial \phi) \, = \, \left\{ (x,y) \in X \times X^{*} \mbox{
: } y \in \partial \phi(x) \right\}$$
for some $\displaystyle \phi \in \Gamma_{0}(X)$. Remark that even in the 
general case of a pair of spaces $X, Y$ in duality, one implication from  this
result is still true, namely if $M$ is cyclically maximal monotone then there 
exists $\displaystyle \phi \in \Gamma_{0}(X)$ such that 
$\displaystyle M  =  Graph(\partial \phi)$. Moreover, by Fenchel inequality 
this means that 
\begin{equation}
b(x,y) = \phi(x) + \phi^{*}(y)
\label{sepab}
\end{equation}
is a bipotential and $M(b) = Graph(\partial \phi)$. Such a bipotential is called
separable.

\begin{rk} Although using the term "graph of $b$" for the set defined by 
(\ref{mb}) seems not adequate, this denomination is used repeatedly in 
previous articles concerning bipotentials, therefore we shall keep it for 
coherence reasons. A motivation for the use of word "graph" comes from 
the relation $M(b) = Graph(\partial \phi)$ in the case of a separable
bipotential.
\label{rem1}
\end{rk}

Bipotentials are related to syncronised convex functions, defined further. 

\begin{defi}
A {\bf sync} (syncronised convex function) is a function $c: X \times Y
\rightarrow 
[0,+\infty]$ with the properties: 
\begin{enumerate}
\item[(a)] for any $x \in X$, if $dom \, c(x, \cdot) \not = \emptyset$ then 
$\displaystyle c(x, \cdot) \in \Gamma_{0}(X)$;  for any $y \in Y$, if $dom \,
c(\cdot, y) \not = \emptyset$ then 
$\displaystyle c(\cdot, y) \in \Gamma_{0}(Y)$; 
\item[(b)] for any $x \in X$, if  $dom \, c(x,\cdot) \not = \emptyset$ and 
the minimum  $\min\left\{ c(x,y) \mbox{ : } y \in
Y \right\}$ exists then this minimum equals $0$; for any $y \in X$, if  $dom \,
c(\cdot, y) \not = \emptyset$ and 
the minimum  $\min\left\{ c(x,y) \mbox{ : } x \in
X \right\}$ exists then this minimum equals $0$. 
\end{enumerate}
\label{defsync}
\end{defi}

\begin{prop}
A function $b: X \times Y \rightarrow
 \bar{\mathbb{R}}$ is a bipotential if and only if the function 
$c: X \times Y \rightarrow
 \bar{\mathbb{R}}$, $c(x,y) = b(x,y) - \langle x, y \rangle$ is a sync. 
\label{psync}
\end{prop}

\paragraph{ Proof.} 
With the notations from the proposition,  
the property (b) definition \ref{def2} for the function $b$   
%%% BEGIN de Saxce
is equivalent with $c: X \times Y  \rightarrow
[0,+\infty]$. The property (a) definition \ref{def2} for the function $b$ 
%%% END de Saxce
 is equivalent with 
property (a) definition \ref{defsync} for the function $c$. Finally, the 
string of equivalences (\ref{equiva})  
from property (c) definition \ref{def2} for the function $b$ 
is equivalent with the following property for the function $c$:
\begin{equation} 
0 \in \partial c(\cdot , y)(x) \ \Longleftrightarrow \ 0 \in \partial c(x,
\cdot)(y)  \ \Longleftrightarrow \ c(x,y) = 0 \ .
\label{equivb}
\end{equation}
But this is just a reformulation of the 
property 
(c) definition \ref{defsync} for the function $b$.  $\quad \square$ 

\begin{rk}
The string of equivalences (\ref{equivb}) justifies the name 
"syncronised convex function", as it expresses the fact that critical 
points of functions $c(x, \cdot)$ are related with critical points of 
functions $c(\cdot, y)$. 
\end{rk}

With the notations from proposition \ref{psync}, we have  
$\displaystyle M(b) = c^{-1}(0)$. Also, for any $x \in X$ and $y \in Y$,
property (a) definition \ref{defsync} of syncs is equivalent with: 
$$epi(c) \cap \left\{x\right\} \times Y \times \mathbb{R} \mbox{ and } 
epi(c) \cap X \times \left\{y\right\}  \times \mathbb{R}$$ 
are closed convex sets, where $epi(c)$ is the epigraph of $c$: 
$$epi(c) = \left\{ (x,y,r) \in X \times Y \times \mathbb{R} \mbox{ : } 
c(x,y) \leq r \right\}$$
For any graph $M \subset X \times Y$, we introduce the sections 
$\displaystyle M (x) \ = \ \left\{ y \in Y \mid (x,y) \in M \right\} \ $ and 
$ M^{*}(y) \ = \ \left\{ x \in X \mid (x,y) \in M \right\} \ $. The following is
definition 3.1 \cite{bipo1}.

\begin{defi}
$M \subset X \times Y$ is a {\bf BB-graph} (bi-convex and 
 bi-closed) if for any $x \in X$ and $y \in Y$ the sections 
$M(x)$ and $M^{*}(y)$ are convex and closed. 
\end{defi}

For any non empty BB-graph $M$ the indicator function $\displaystyle \chi_{M}$ is 
obviously a sync.

\section{Existence and non uniqueness of the bipotential}

Let a constitutive law be given by a graph $M$. Does it admit a bipotential? 
The existence problem is easily settled by the following result (theorem 3.2 
\cite{bipo1}).

\begin{thm}
 Given  a non empty set $M \subset X \times Y$, there is a bipotential $b$ 
such that $M=M(b)$ if and only if $M$ is a BB-graph. 
 \label{thm1}
 \end{thm}

 To any  BB-graph $M$ is  associated  the sync $\displaystyle 
 \chi_{M}$. To this sync corresponds the bipotential  
 \begin{equation}
  b_{\infty} (x,y) = \left\langle x,y \right\rangle  + \chi_M (x,y) .
  \label{binf}
  \end{equation} 
In particular, this shows that to a BB-graph we may associate more than one 
bipotential. Indeed, if $M$ is maximal cyclically monotone  then there are 
two different bipotentials associated to it. First there is the separable 
bipotential (\ref{sepab}), thus $M = M(b) = \, Graph(\partial \phi)$. 
But if $M$ is a BB-graph then  $\displaystyle \chi_{M}$ is a sync, 
which implies that $\displaystyle M = M(b_{\infty})$ with  $\displaystyle
b_{\infty}$ the bipotential defined at 
(\ref{binf}). Therefore maximal cyclically monotone 
graphs admit at least two distinct bipotentials.

The graph  alone is not  sufficient to uniquely define 
the bipotential. However, remark that syncs expressed as indicator functions
don't seem useful as constitutive laws, as they are somehow trivial. 
Therefore we would like to be able to
construct more interesting bipotentials, for example we want a method 
of construction of bipotentials which associates to a maximal cyclically
monotone set a separable bipotential. 

Nevertheless the trivial indicator functions will prove to be useful 
in connection to blurred constitutive laws.

\section{Bipotential convex covers}

 Theorem \ref{thm1} does not give a satisfying bipotential for a given multivalued constitutive law, because the bipotential $b_{\infty}$ is somehow degenerate. We would like to be able to find a bipotential $b$ which is not everywhere infinite outside the graph 
$M$. We saw that the graph alone is not sufficient to construct interesting 
bipotentials. We need more information to start from. This is provided by the notion of bipotential convex cover. 

Let $Bp(X,Y)$ be the set of all bipotentials $b: X \times Y \rightarrow
\bar{\mathbb{R}}$. We shall need the following definitions.
\begin{defi}
Let $\Lambda$ be an arbitrary non empty set and $V$ a real vector space. The 
function $f:\Lambda\times V \rightarrow \bar{\mathbb{R}}$ is 
{\bf implicitly  convex} if for any two elements 
$\displaystyle (\lambda_{1}, z_{1}) , 
(\lambda_{2},  z_{2}) \in \Lambda \times V$ and for any two numbers 
$\alpha, \beta \in [0,1]$ with $\alpha + \beta = 1$ there exists 
$\lambda  \in \Lambda$ such that 
\begin{equation}
   f(\lambda, \alpha z_{1} + \beta z_{2}) \ \leq \ \alpha 
   f(\lambda_{1}, z_{1}) + \beta f(\lambda_{2}, z_{2}) \quad .
\label{ineqimpl} 
\end{equation} 
\label{defimpl}
\end{defi}

\begin{defi}   A {\bf bipotential
convex cover} of the non empty set $M$ is a function   
$\displaystyle \lambda \in \Lambda \mapsto b_{\lambda}$ from  $\Lambda$ with 
values in the set  $Bp(X,Y)$, with the 
properties:
\begin{enumerate}
\item[(a)] The set $\Lambda$ is a non empty compact topological space, 
\item[(b)] Let $f: \Lambda \times X \times Y \rightarrow \mathbb{R} \cup
\left\{ + \infty \right\}$ be the function defined by 

$$f(\lambda, x, y) \ = \ b_{\lambda}(x,y) .$$

Then for any $x \in X$ and for any $y \in Y$ the functions 
$f(\cdot, x, \cdot): \Lambda \times Y \rightarrow \bar{\mathbb{R}}$ and 
$f(\cdot, \cdot , y): \Lambda \times X \rightarrow \bar{\mathbb{R}}$ are  lower 
semi continuous  on the product spaces   $\Lambda \times Y$ and respectively 
$\Lambda \times X$ endowed with the standard topology, 
\item[(c)] We have $\displaystyle M  \ = \  \bigcup_{\lambda \in \Lambda} 
M(b_{\lambda})$.
\item[(d)] with the notations from point (b), the functions $f(\cdot, x, \cdot)$ 
and $f(\cdot, \cdot , y)$ are implicitly convex in the sense of Definition 
\ref{defimpl}.
\end{enumerate}
\label{defcover}
\end{defi}

A bipotential convex cover is in some sense described by the collection 
$\displaystyle \left\{ b_{\lambda} \mbox{ : } \lambda \in \Lambda \right\}$. 
This is this point of view that we will adopt in the sequel. The next result 
defines under which conditions the notion of bipotential convex cover is 
independent of the choice of the parameterization 
\cite{bipo2}. 

\begin{prop}
Let $\displaystyle \lambda \in \Lambda \mapsto b_{\lambda} \in Bp(X,Y)$ be a
bipotential convex cover and $g: \Lambda \rightarrow \Lambda$ be a homeomorphism. Then 
 $\displaystyle \lambda \in \Lambda \mapsto b_{g(\lambda)} \in Bp(X,Y)$ is a
 bipotential convex cover. 
\end{prop}

The next theorem, (4.6 \cite{bipo3}), is the key result needed further.

\begin{thm} Let $\displaystyle \lambda \mapsto b_{\lambda}$ be a bipotential 
convex cover  of 
 the graph $M$ and $b: X \times Y \rightarrow \bar{\mathbb{R}}$ defined by
\begin{equation}
b(x,y) \ = \ \inf \left\{ b_{\lambda}(x,y) \ \mid \  \lambda \in \Lambda \right\} \ . 
\end{equation}
Then $b$ is a bipotential and $M=M(b)$. 
\label{thm2}
\end{thm}

An inferior envelope of  convex functions is not generally convex.  
The property (d)  of the Definition \ref{defcover} is essential to ensure 
the convexity properties of $b$.

\section{Blurred BB-graphs}

In many practical situations, indeterminacies affect the mechanical behaviour. 
In other words, we tolerate indeterminacy of the constitutive law. 

We shall represent the indeterminacy by a set 
$A \subset X \times Y, (0,0) \in A$ and we
shall suppose that it is a BB-graph. This hypothesis is justified by the
following examples. 

Suppose $X$ is a reflexive Banach space  and 
$\langle x, y\rangle = y(x)$. We shall denote by $\|\cdot\|$ both norms, 
in $X$ and in $Y$.  Let 
$\displaystyle A = \left\{0\right\} \times \bar{B}_{Y}(\varepsilon) = \left\{ (0,y) \mbox{ : } 
\|y \| \leq \varepsilon \right\}$. 
This set is a BB-graph and represents the indeterminacy $\varepsilon$ in the
norm of $y$, for given $x$. 

In the same setting we can see  $X \times Y$ as a  
normed vector space with the norm 
$$\|(x,y)\| = \left( \|x\|^{p} + \|y\|^{p}\right)^{1/p}$$ 
for a $p \geq 1$. Take then 
$\displaystyle A = \left\{ (x,y) \in X \times Y \mbox{ : } \|(x,y)\| \leq
\varepsilon \right\}$. 
This set, which is a BB-graph, represents the indeterminacy $\varepsilon$ in the
norm of the pair $(x,y)$. 

Let $c: X \times Y \rightarrow [0,+\infty]$ be a sync and $\displaystyle 
M = c^{-1}(0)$ the graph of the associated bipotential $b$. Suppose that 
the graph $M$ represents a constitutive law. The constitutive law with 
indeterminacy represented by a BB-graph $A \subset X \times Y$ has then the
graph 
$$M + A = \left\{ (x,y) \in X \times Y \mbox{ : } x = x'+x", y = y'+y", (x', y')
\in M , (x",y") \in A \right\}$$
For example, if $\displaystyle c(x,y) = \phi(x) + \phi^{*}(y) - \langle x, y
\rangle$, with $\displaystyle \phi \in \Gamma_{0}(X)$ then $M$ is the graph 
of the subdifferential $\partial \phi$. If we take the indeterminacy 
$\displaystyle A = \left\{0\right\} \times \bar{B}_{Y}(\varepsilon)$ then 
\begin{equation}
M + A = \left\{ (x,y) \in X \times Y \mbox{ : } \exists a \in Y, \|a\| \leq
\varepsilon, \, y + a \in \partial \phi(x) \right\}
\label{neednewc}
\end{equation}
A natural function associated to $M+A$ is the inf-convolution 
$$c_{A}(x,y) = \left(c \nabla \chi_{A} \right)(x,y) = \inf \left\{ 
c(x',y') + \chi_{A}(x",y") \mbox{ : } x'+x"=x, y'+y"=y\right\}$$
By the properties of the inf-convolution we have 
$$epi(c_{A}) = epi(c) + epi(\chi_{A}) = epi(c) + \left( A \times (0,
+\infty)\right)$$
By definition of $\displaystyle c_{A}$ we have 
 $\displaystyle c^{-1}(0) + A = M+A \subset c_{A}^{-1}(0)$. Under supplementary 
 hypothesis on $A$ and $c$ we have the equality  
  $\displaystyle  M+A =  c_{A}^{-1}(0)$. Let us present such hypothesis. 
 
 On the space $X \times Y$ we consider the convergence $\displaystyle 
 (x_{h}, y_{h}) \rightarrow (x,y)$ (as $h \rightarrow \infty$) if 
 $\displaystyle x_{h}$ converges weakly to $x$ and $\displaystyle y_{h}$
 converges weakly* to $y$. If  for any $(x,y) \in X \times Y$ 
 the function $(x", y") \in A \mapsto c(x - x", y - y")$ is lsc with respect to
 this convergence  and  
 $A$ is compact with respect to this convergence then $\displaystyle  M+A = c_{A}^{-1}(0)$. Indeed, for any 
 $h \in \mathbb{N}^{*}$ there is $\displaystyle (x_{h}, y_{h}) \in A$ such that 
 $\displaystyle 0 \leq c(x - x_{h}, y - y_{h}) \leq 1/h$. Because $A$ is
 compact, we can extract a subsequence converging to a pair $(x', y') \in A$; 
 by the lsc of $c$ we find that 
 $\displaystyle 0 \leq c(x - x', y - y') \leq \liminf_{h \rightarrow \infty} 
 c(x - x_{h}, y - y_{h}) = 0$, therefore $\displaystyle c(x - x', y - y') = 0$. 
 This proves the inclusion $\displaystyle  c_{A}^{-1}(0) \subset M + A$,
 therefore we get the desired equality of sets. 
 
 A particular case when this hypothesis is obviously 
 true is  the one mentioned in relation (\ref{neednewc}), 
 presented further in  section \ref{sec6}, that is 
$\displaystyle A = \left\{ 0 \right\} \times \bar{B}_{Y}(\varepsilon)$ and 
$\displaystyle c(x,y) = \ \phi(x) + \phi^{*}(y) - \langle x, y \rangle$, with 
$\displaystyle \phi \in \Gamma_{0}(X)$. Another
situation when the mentioned hypothesis is true is the one when 
$X$ and $Y$ are finite dimensional and $A$ is bounded.

\begin{defi}
Consider a BB-graph  $M \subset X \times Y$ and another BB-graph $A \subset X \times
Y$ such that $(0,0) \in A$. We say that $M$ admits the blurring $A$ 
if   $M+A$  is a BB-graph. 

Let $c:  X \times Y \rightarrow [0,+\infty]$ be a sync and $A \subset X \times
Y$ be a BB-graph such that $(0,0) \in A$. We say that $c$ admits the blurring
$A$ if $\displaystyle c_{A} = c \nabla \chi_{A}$ is a sync and 
$\displaystyle c_{A}^{-1}(0) = c^{-1}(0) + A$. 
\label{ablurred}
\end{defi}
 
\paragraph{Example 1.} (Blurred elasticity.) We  take $X=Y=\mathbb{R}^{n}$, 
the duality product is the usual scalar product in $\mathbb{R}^{n}$ and
$\|\cdot\|$ is the usual norm. 
%%% BEGIN de Saxce
We consider the elastic linear law $y = K x$ which is 
%%% END de Saxce
the most simple example of linear elastic law where the dual variables 
$x$ and $y$ are vectors, and the "elastic modulus" $K > 0$. To this law is associated the graph: 
$$M = \left\{ (x,y) \in X \times Y \mbox{ : } y = K x \right\}$$
This graph is maximal cyclically monotone and it admits the sync: 
$$c(x,y) = \frac{K}{2} \|x\|^{2} + \frac{1}{2K}\|y\|^{2} - \langle x,y \rangle$$
Let $\varepsilon > 0$ and 
$\displaystyle A = \left\{0\right\} \times \left\{ y \in Y \mbox{ : } \|y\| \leq \varepsilon
\right\}$. 
Then we have: 
$$M + A = \left\{ (x,y) \in X \times Y \mbox{ : } \|y - Kx\| \leq \varepsilon
\right\}$$
which is a BB-graph, therefore $M$ admits the blurring $A$. Moreover, after some
computations  we get that  
$$c_{A}(x,y) = \left( c\nabla \chi_{A} \right)(x,y) = \frac{1}{2K} \left(
\left(\|y-Kx\| - \varepsilon\right)_{+}\right)^{2}$$
with the notation $\displaystyle z \in \mathbb{R} \mapsto \left(z\right)_{+} = \max(z,0)$.  It is easy to check that $\displaystyle c_{A}$ is a sync. 
Therefore the graph $M$  and the sync $c$ admit the blurring $A$. A similar
computation can be done in the case $K$ is a strictly positive definite matrix, 
only that the expression of $\displaystyle c_{A}$ is more complex. 

\paragraph{Example 2.} (A BB-graph made of two points.) In the setting of the
previous example, consider this time $\displaystyle 
M = \left\{ (x_{1}, y_{1}), (x_{2}, y_{2}) \right\}$, with $\displaystyle 
x_{1} \not = x_{2}$ and $y_{1} \not = y_{2}$. This is a BB-graph 
(although not a very interesting one). We take $A$ as previously. Remark that 
if $\displaystyle \|y_{1} - y_{2} \| \leq 2 \varepsilon$ then $M+A$ is not 
bi-convex. In this case $M$ does not admit the blurring $A$. 

\section{Bipotentials for blurred maximal cyclically monotone sets}
\label{sec6}

In this section $X$ is a reflexive Banach space, $Y$ the dual space and 
the symbol $\|\cdot\|$ is used to denote the norm in $X$ or the dual norm 
in $Y$.

As the previous example looks somehow degenerate, we
might hope that at least in the case $M$ is a maximal cyclically monotone set 
then $M$ admits the blurring $A$, where $A$ is defined as in the examples above.

The following proposition gives a necessary and sufficient condition for this. 

\begin{prop}
Let  $\varepsilon > 0$,  
 $M = Graph(\partial \phi)$  be a maximal cyclically monotone set with  
$\displaystyle \phi \in \Gamma_{0}(X)$, and $A \subset X \times Y$ be defined 
by $\displaystyle A = \left\{ 0 \right\} \times \bar{B}_{Y}(\varepsilon)$. 
The set $M$ admits the blurring $A$ if and only if the following condition is true: 
\begin{equation}
\mbox{ for any } y \in Y \quad  \mbox{ the set } \bigcup_{\| \bar{y}-y\|\leq \varepsilon} \partial 
\phi^{*}(\bar{y}) \quad \mbox{ is convex} 
\label{newc}
\end{equation}
\label{pnewc}
\end{prop}

\paragraph{Proof.}
The expression of $M+A$ is given in (\ref{neednewc}). The set $M+A$ is a
BB-graph if and only if for any $y \in Y$ the set $\displaystyle \left\{ x \in X \mbox{ : } 
(x,y) \in M+A \right\}$ is convex and closed. A simple computation shows that: 
\begin{equation} 
\left\{ x \in X \mbox{ : } 
(x,y) \in M+A \right\} = \bigcup_{\| \bar{y}-y\|\leq \varepsilon} \partial 
\phi^{*}(\bar{y})
\label{weneed2}
\end{equation} 
This set is closed, for any $\displaystyle \phi \in \Gamma_{0}(X)$ and 
$\varepsilon > 0$. Indeed, let $\displaystyle x_{h}$, $\displaystyle h \in
\mathbb{N}$ be a sequence in this set, converging to $x \in X$ as $h$ goes to 
infinity. For any $\displaystyle h \in \mathbb{N}$ we have $\displaystyle 
y +  a_{h}  \in \partial \phi( x_{h})$, with $\displaystyle \|a_{h}\| \leq
\varepsilon$. Equivalently, for  any $\displaystyle h \in \mathbb{N}$ there
exists $\displaystyle a_{h} \in Y$, $\displaystyle \|a_{h}\| \leq
\varepsilon$ such that $\displaystyle 
\phi( x_{h}) < + \infty$ and for any $x' \in X$ we have 
\begin{equation}
\phi(x') - \langle x', y + a_{h} \rangle \geq \phi(x_{h}) -  
\langle x_{h},  y +  a_{h}\rangle
\label{weneed}
\end{equation}
We may suppose without reduction of generality that 
$\displaystyle a_{h}$ weakly* converges to $a \in Y$, with $\|a \| \leq
\varepsilon$. For a fixed, but arbitrary $x' \in X$ we pass to the limit in 
(\ref{weneed}) and we get: 
$$\phi(x')  - \langle x', y + a \rangle \, = \, \lim_{h \rightarrow \infty} 
\left[ \phi(x') - \langle x', y + a_{h} \rangle \right] \, \geq \, $$ 
$$\geq \, \liminf_{h \rightarrow \infty} \left[ \phi(x_{h}) -  
\langle x_{h},  y +  a_{h}\rangle \right] \, \geq \, 
\phi(x) -  
\langle x,  y +  a\rangle$$ 
We proved that $\displaystyle y + a \in \partial \phi(x)$, with  $\|a \| \leq
\varepsilon$.

The relation (\ref{weneed2}) lead us to the conclusion that 
 condition (\ref{newc}) is necessary and sufficient for the set 
$M$ to be convex. $\quad \square$

\paragraph{Example 3.} 
The condition (\ref{newc}) is not true for any $\displaystyle 
\phi \in \Gamma_{0}(X)$, at least when $X=Y=\mathbb{R}^{n}$ and $n \geq 2$.  
Indeed, for $n=2$  let us  take $\displaystyle \phi^{*}(y) = \chi_{F}(y)$ where $F \subset Y$ 
is the closed convex cone
$$F = \left\{ (y_{1}, y_{2}) \in \mathbb{R}^{2} \mbox{ : } \mid y_{2}\mid \leq 
\alpha y_{1} \right\}$$
with $\alpha \in (0,1)$. Obviously then $\displaystyle 2 \alpha < 1+
\alpha^{2}$. Remark that the boundary of $F$ is made by two half lines: 
$\displaystyle h_{1}, h_{2}$ of equation  
$\displaystyle  \pm y_{2} = \alpha y_{1}$, $\displaystyle y_{1} \geq 0$. These
half lines have normals denoted by $\displaystyle n_{1}, n_{2}$, which are not 
proportional one with another. 

Let $M = Graph(\partial \phi)$. We shall choose now $\varepsilon$ such that there
exists $\displaystyle y_{1} > 0$ with 
$$\frac{2\alpha}{\sqrt{1+\alpha^{2}}} y_{1} < \varepsilon < y_{1}
\sqrt{1+\alpha^{2}}$$ 
Let $\displaystyle y_{2} = \alpha y_{1}$. We have then $\displaystyle 
(y_{1}, y_{2}) \in F$ and 
$$ (y_{1}, y_{2}) + \bar{B}_{Y}(\varepsilon) \cap h_{1} \not = \emptyset \quad
, \, (y_{1}, y_{2}) + \bar{B}_{Y}(\varepsilon) \cap h_{2} \not = \emptyset$$
and $\displaystyle (0,0) \not \in (y_{1}, y_{2}) + \bar{B}_{Y}(\varepsilon)$.
Then 
$$\bigcup_{\| \bar{y}-y\|\leq \varepsilon} \partial 
\phi^{*}(\bar{y}) = \left\{ \lambda n_{1}, \lambda n_{2} \mbox{ : } \lambda \geq
0 \right\}$$ 
which is not a convex set. Therefore $M + A$ is not a BB-graph. 

In the following consider an arbitrary  $\displaystyle \phi \in 
\Gamma_{0}(X)$ and $M = Graph(\partial \phi)$.  To 
 $\phi$ is associated the separable bipotential $\displaystyle 
 b(x,y) = \phi(x) + \phi^{*}(y)$ and the sync $c(x,y) = b(x,y) - \langle x, y
 \rangle$. We have then $\displaystyle M = M(b) = c^{-1}(0)$. 
 
 Take  $A \subset X \times Y$, 
 $\displaystyle A = \left\{ 0 \right\} \times \bar{B}_{Y}(\varepsilon)$ and 
 let $\displaystyle c_{A}(x,y) = \left(c \nabla \chi_{A} \right)(x,y)$. We saw
 that $\displaystyle c_{A}^{-1}(0) = M+A$. We want
 to know if $\displaystyle c_{A}$ is a sync, equivalently if 
 $$b_{A}(x,y) = \langle x, y \rangle + \inf \left\{ 
c(x',y') + \chi_{A}(x",y") \mbox{ : } x'+x"=x, y'+y"=y\right\}$$
is a bipotential. Proposition \ref{pnewc} gives us the necessary condition 
(\ref{newc})  
for this to be true, because if $M+A$ is not a BB-graph then $\displaystyle 
c_{A}$ cannot be a sync. 

Further we shall give a necessary and sufficient condition for 
$\displaystyle c_{A}$ to be a sync. Before this, we want to comment on the 
form of the function $\displaystyle b_{A}$. We start from the remark that for 
any $a \in Y$ and any sync $c: X \times Y \rightarrow [0,+\infty]$ the 
function defined by $\displaystyle c_{a}(x,y) = c(x, y-a)$ is also a sync. 
In particular,
for any $\displaystyle \phi \in \Gamma_{0}(X)$ the function
$$b_{a}(x,y) = \phi(x) + \phi^{*}(y-a) + \langle x, a\rangle$$
is a (separable) bipotential.

We see that the function $\displaystyle b_{A}$ has the expression 
$$b_{A}(x,y) = \inf \left\{ \phi(x) + \phi^{*}(y-a) + \langle x, a\rangle \mbox{
: } \|a\| \leq \varepsilon \right\}$$
therefore, $\displaystyle b_{A}$ is an infimum of bipotentials: 
$$b_{A}(x,y) = \inf \left\{ b_{a}(x,y) \mbox{ : } \|a\| \leq \varepsilon
\right\}$$
If the function $\displaystyle a \in \bar{B}_{Y}(\varepsilon) \mapsto b_{a}$ 
is a bipotential convex cover then $\displaystyle b_{A}$ is a bipotential, by
theorem \ref{thm2}.

\begin{thm}
 With the previous notations $\displaystyle b_{A}$ is a bipotential such that 
 $\displaystyle M(b_{A}) = M(\phi, \varepsilon)$ if and 
only if  for any $y \in Y$ the function $\displaystyle f(\cdot, \cdot, y): 
\bar{B}_{Y}(\varepsilon) \times X 
\rightarrow \bar{\mathbb{R}}$ defined by 
$$f(a, x, y) = \phi(x) + \phi^{*}(y - a) + \langle x, a \rangle$$ 
is implicitly convex. 
\label{maithm}
\end{thm}

\paragraph{Proof.}
If $\displaystyle b_{A}$ is a bipotential then for any $ y \in Y$ the function 
$\displaystyle b_{A}(\cdot, y)$ is convex. By reasoning as in the proof of 
corollary \ref{coro2}, it follows that $\displaystyle f(\cdot, \cdot, y)$ is
implicitly convex. 

We have to prove now the inverse implication.  We
shall prove that the function $\displaystyle a \in \bar{B}_{Y}(\varepsilon) \mapsto b_{a}$ 
is a bipotential convex cover of the set $M(\phi, \varepsilon)$. 

Let us denote
 $\displaystyle \Lambda = \bar{B}_{Y}(\varepsilon)$. This set, 
endowed with  the weak  topology from $Y$ is a compact topological space. Then
the function $\displaystyle f:  \Lambda \times X 
\rightarrow \bar{\mathbb{R}}$ defined by 
$$f(a, x, y) = \phi(x) + \phi^{*}(y - a) + \langle x, a \rangle$$ 
has the following properties. For any $x \in X$ the function 
$f(\cdot, x, \cdot)$ is convex and lower semicontinuous, with respect to the 
product topology (weak topology on $\Lambda$ times the strong topology on $Y$). 
Also, for any $y \in Y$ the function $f(\cdot, \cdot, y)$ is implicitly convex 
by hypothesis, and also lower semicontinous with respect to the 
product topology. Indeed, this is true because: $\phi$ is strongly lower
semicontinuous, the function  $a \in \Lambda \mapsto \phi^{*}(y-a)$ is weakly 
lower semicontinuous and for any sequence $\displaystyle (a_{h}, x_{h}) \in \Lambda \times X$ such that $\displaystyle 
a_{h}$ weakly converges to $ a \in \Lambda$ and $\displaystyle x_{h}$ strongly
converges to $x \in X$, then $\displaystyle \langle x_{h}, a_{h} \rangle$
converges to $\langle x, a \rangle$. 

It follows that for any $(x,y) \in X \times Y$ there exists $\bar{a} \in \Lambda$ such
that 
$$b_{A}(x,y) = \min_{a \in \Lambda} f(a, x, y)$$
We deduce that $\displaystyle M(\phi, \varepsilon) = \bigcup_{a \in \lambda}
M(b_{a})$. We proved that $\displaystyle a \in \bar{B}_{Y}(\varepsilon) \mapsto b_{a}$ 
is a bipotential convex cover of the set $M(\phi, \varepsilon)$. 
By theorem \ref{thm2} we obtain the conclusion. 
 $\quad \square$

\vspace{.5cm}

\paragraph{Final remark.} Suppose $X = \mathbb{R}$. Then condition 
 (\ref{newc}) should take a simpler form, because in one dimension 
 convex is the same as connected and (\ref{newc}) is just a kind of 
 Darboux property for subgradients. Is it true for any  
 $\phi \in \Gamma_{0}(X)$? 

\vspace{\baselineskip}

\end{document}